\documentclass[11pt]{article}

\usepackage[a4paper,margin=1in]{geometry}
\usepackage{amsmath,amssymb,amsthm,mathtools,bm}
\usepackage{enumitem}
\usepackage{hyperref}
\usepackage{mathrsfs}
\usepackage{array}
\usepackage{booktabs}
\usepackage{xcolor}

\usepackage{authblk} 
\usepackage[numbers,sort&compress]{natbib} 

\usepackage{tikz}
\usetikzlibrary{
  arrows.meta,
  calc,
  positioning,
  decorations.pathreplacing,
  patterns,
  fit,
  backgrounds,
  intersections,
  shapes.geometric  
}

\hypersetup{
  colorlinks=true,
  linkcolor=blue!60!black,
  citecolor=blue!60!black,
  urlcolor=blue!60!black
}

\newtheorem{theorem}{Theorem}[section]

\theoremstyle{definition}
\newtheorem{definition}[theorem]{Definition}
\newtheorem{remark}[theorem]{Remark}

\newcommand{\R}{\mathbb{R}}

\newcommand{\T}{\mathsf{T}}
\newcommand{\Gr}{\operatorname{Gr}}

\newcommand{\Lag}{\mathcal{L}}
\newcommand{\MPEC}{\mathrm{MPEC}}
\newcommand{\NLP}{\mathrm{NLP}}

\newcommand{\C}{\mathcal{C}}
\newcommand{\F}{\mathcal{F}}

\tikzset{
  axis/.style={->, thick},
  ray/.style={-Latex, thick},
  mybox/.style={draw, rounded corners, thick, inner sep=6pt, align=center},
  myellipse/.style={draw, ellipse, thick, align=center, inner sep=4pt},
  point/.style={circle, fill=black, inner sep=1.6pt},
  critical/.style={very thick, dashed},
  feasible/.style={thick},
  labelbox/.style={draw, rounded corners, fill=gray!8, inner sep=4pt, align=left}
}

\title{Optimization Workshop Notes for Mathematical Programming with Equilibrium Constraints (MPECs): \\ Second-Order Optimality Conditions}
\author{Jiguang Yu\thanks{Email: jyu678@bu.edu}}
\affil{College of Engineering, Boston University, Boston, 02215, MA, USA}

\begin{document}
\maketitle

\begin{abstract}
In this workshop, we present a compact but rigorous introduction to second-order
optimality conditions for mathematical programs with equilibrium constraints
(MPECs). We start from the classical nonlinear programming template,
then explain why it fails in the equilibrium-constrained setting, and develop
the three main viewpoints used in the literature:
(i) multiplier-based conditions,
(ii) implicit-programming conditions based on the solution map of the lower-level equilibrium system,
and (iii) piecewise-programming conditions obtained by decomposing complementarity structure into smooth pieces.
The emphasis is on conceptual structure, critical cones, strong regularity,
and the exact role of curvature terms.
\end{abstract}


\section{Introduction}

An MPEC is an optimization problem whose feasible set is constrained by
an equilibrium system, often a variational inequality (VI), a nonlinear
complementarity problem (NCP), or a KKT system. A prototypical form is
\begin{equation}\label{eq:mpec-general}
\begin{aligned}
\min_{x,y} \quad & f(x,y)\\
\text{s.t.}\quad & (x,y)\in Z,\\
& y\in S(x),
\end{aligned}
\end{equation}
where $S(x)$ is the solution set of an equilibrium problem parameterized by $x$.

The key analytic difficulty is that even when the data are smooth, the graph
$\Gr S$ is typically only \emph{piecewise smooth} or \emph{piecewise polyhedral},
and standard nonlinear-programming constraint qualifications often fail.
As a consequence, second-order theory for MPECs cannot simply be imported
verbatim from smooth NLP.

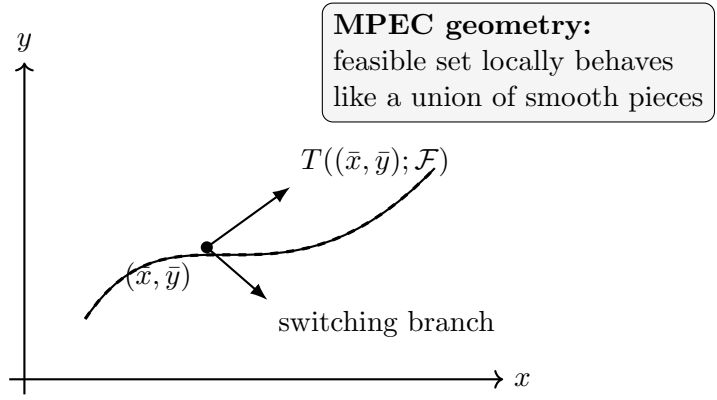
\begin{figure}[ht]
\centering
\begin{tikzpicture}[scale=1.0]
  \draw[axis] (-0.2,0) -- (6.3,0) node[right] {$x$};
  \draw[axis] (0,-0.2) -- (0,4.2) node[above] {$y$};

  \draw[feasible] (0.8,0.8) .. controls (2,2.6) and (3.3,0.6) .. (5.4,2.8);
  \draw[critical] (0.8,0.8) .. controls (2,2.6) and (3.3,0.6) .. (5.4,2.8);

  \node[point, label=below left:{$(\bar x,\bar y)$}] at (2.4,1.75) {};
  \draw[ray] (2.4,1.75) -- (3.5,2.55) node[above right] {$T((\bar x,\bar y);\F)$};
  \draw[ray] (2.4,1.75) -- (3.2,1.05) node[below right] {switching branch};

  \node[labelbox] at (6.5,4.2) {
    \textbf{MPEC geometry:}\\
    feasible set locally behaves\\
    like a union of smooth pieces
  };
\end{tikzpicture}
\caption{Local picture: the feasible set of an MPEC is typically branchwise smooth rather than globally smooth.}
\end{figure}

The notes focus on the following themes:

\begin{enumerate}[label=(\roman*)]
\item the role of the \emph{critical cone} as the carrier of second-order information;
\item why complementarity constraints are better viewed as \emph{disjunctive} than as ordinary smooth equalities;
\item how \emph{strong regularity} of the lower-level equilibrium system enables an implicit-function approach;
\item how multiplier-based and implicit-programming conditions complement one another.
\end{enumerate}

\section{Classical benchmark: second-order theory for NLP}

Consider the nonlinear program
\begin{equation}\label{eq:nlp}
\begin{aligned}
\min_{x\in\R^n}\quad & f(x)\\
\text{s.t.}\quad & g_i(x)=0,\quad i=1,\dots,\ell,\\
& g_i(x)\le 0,\quad i=\ell+1,\dots,\ell+m,
\end{aligned}
\end{equation}
with $f$ and $g_i$ of class $C^2$.

\subsection{Stationarity and critical cone}

Let $\bar x$ be feasible. The active set is
\[
A(\bar x)
:=
\{1,\dots,\ell\}
\cup
\{i>\ell : g_i(\bar x)=0\}.
\]
If $\bar x$ is stationary, there exists a multiplier vector
$\bar\pi\in\R^\ell\times\R^m_+$ such that
\[
\nabla f(\bar x)+\sum_{i=1}^{\ell+m}\bar\pi_i \nabla g_i(\bar x)=0,
\qquad
\bar\pi_i g_i(\bar x)=0\;\; (i>\ell).
\]

\begin{definition}[Critical cone]
Assume $\bar x$ is stationary. The critical cone is
\[
\C(\bar x)
:=
\left\{
d\in T(\bar x;\F):
\nabla f(\bar x)^\T d = 0
\right\},
\]
where $T(\bar x;\F)$ is the tangent cone to the feasible set $\F$ of \eqref{eq:nlp}.
\end{definition}

\begin{figure}[ht]
\centering
\begin{tikzpicture}[scale=1.2, >=Stealth]
  \tikzset{
    axis/.style={->, gray, thin},
    point/.style={circle, fill=black, inner sep=1.5pt},
    feasible/.style={thick, fill=gray!10},
    ray/.style={->, blue!70!black, thick},
    critical/.style={->, red!70!black, very thick},
    labelbox/.style={draw, fill=white, font=\footnotesize, align=center, rounded corners=2pt}
  }

  \draw[axis] (-0.5,0) -- (6.0,0) node[right] {$x_1$};
  \draw[axis] (0,-0.5) -- (0,5.0) node[above] {$x_2$};

  \coordinate (Xbar) at (1.0,1.0);
  \fill[gray!10] (Xbar) -- (5.0,1.0) -- (3.0,4.5) -- cycle;
  \draw[thick] (5.0,1.0) -- (Xbar) -- (3.0,4.5);
  
  \node[point, label=below left:{$\bar{x}$}] at (Xbar) {};

  \draw[ray] (Xbar) -- (4.5,1.0) node[below, pos=0.8, font=\scriptsize] {$\nabla g_1(\bar x)^\T d \le 0$};
  \draw[ray] (Xbar) -- (2.6,3.8) node[left, pos=0.8, rotate=55, font=\scriptsize] {$\nabla g_2(\bar x)^\T d \le 0$};

  \draw[->, ultra thick, green!50!black] (Xbar) -- (0.2,2.5) node[above] {$\nabla f(\bar x)$};

  \draw[critical] (Xbar) -- (3.5,2.2) node[right] {$\mathcal{C}(\bar x)$};
  \node[font=\tiny, gray] at (2.2,1.8) {$\nabla f(\bar x)^\T d = 0$};

  \node[labelbox] at (5,5) {
    \textbf{Classical NLP}\\
    Critical directions satisfy:\\
    1. Feasibility ($d \in \mathcal{T}$)\\
    2. Stationarity ($\nabla f^\T d = 0$)
  };

\end{tikzpicture}
\caption{Tangent cone and critical cone at a boundary stationary point.}
\end{figure}
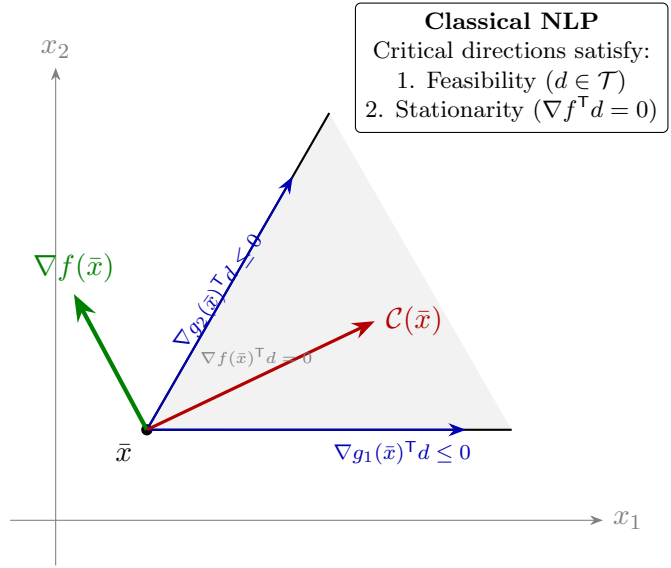

Under MFCQ, one has the familiar linearized representation
\[
\C(\bar x)
=
\left\{
d:
\nabla g_i(\bar x)^\T d=0\;\; (i=1,\dots,\ell),\;\;
\nabla g_i(\bar x)^\T d\le 0\;\; (i\in A(\bar x)\setminus\{1,\dots,\ell\}),\;\;
\nabla f(\bar x)^\T d=0
\right\}.
\]

\subsection{Lagrangian and second-order conditions}

Define the NLP Lagrangian
\[
\Lag_{\NLP}(x,\pi)
=
f(x)+\sum_{i=1}^{\ell+m}\pi_i g_i(x).
\]

The classical template is:

\begin{theorem}[Classical second-order conditions]
Suppose $\bar x$ is a stationary point of \eqref{eq:nlp} and MFCQ holds at $\bar x$.
\begin{enumerate}[label=(\roman*)]
\item If $\bar x$ is a local minimizer, then for every $d\in\C(\bar x)$,
\[
\max_{\pi\in M(\bar x)} d^\T \nabla^2_{xx}\Lag_{\NLP}(\bar x,\pi)\, d \ge 0.
\]
\item If for every nonzero $d\in\C(\bar x)$,
\[
\max_{\pi\in M(\bar x)} d^\T \nabla^2_{xx}\Lag_{\NLP}(\bar x,\pi)\, d > 0,
\]
then $\bar x$ is a strict local minimizer.
\end{enumerate}
\end{theorem}

\section{Why MPEC is different}

A simple complementarity relation
\[
0\le y \perp \lambda \ge 0
\quad\Longleftrightarrow\quad
y\ge 0,\;\lambda\ge 0,\; y_i\lambda_i=0\;\; \forall i
\]
looks smooth if written as the bilinear equations $y_i\lambda_i=0$.
But this representation is deceptive.

\begin{remark}[Disjunctive interpretation]
Each pair $(y_i,\lambda_i)$ satisfies
\[
y_i=0 \quad \text{or} \quad \lambda_i=0.
\]
Thus the feasible set is more naturally seen as the union of finitely many smooth pieces.
This piecewise structure explains why the second derivatives of the bilinear
equalities $y_i\lambda_i=0$ do \emph{not} play the same role they would in an ordinary NLP.
\end{remark}

Two consequences follow immediately:

\begin{enumerate}[label=(\roman*)]
\item naive KKT theory for the full smooth reformulation is unreliable because MFCQ typically fails;
\item second-order conditions must be derived either by exploiting the piecewise geometry directly,
or by using regularity of the lower-level solution map.
\end{enumerate}

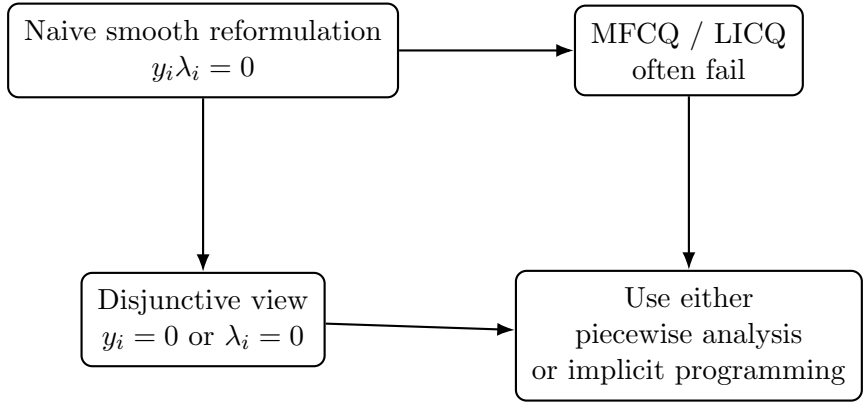
\begin{figure}[ht]
\centering
\begin{tikzpicture}[node distance=2.3cm]
  \node[mybox] (smooth) {Naive smooth reformulation\\$y_i\lambda_i=0$};
  \node[mybox, right=of smooth] (fails) {MFCQ / LICQ\\often fail};
  \node[mybox, below=of smooth] (disj) {Disjunctive view\\$y_i=0$ or $\lambda_i=0$};
  \node[mybox, below=of fails] (tools) {Use either\\piecewise analysis\\or implicit programming};

  \draw[-Latex, thick] (smooth) -- (fails);
  \draw[-Latex, thick] (smooth) -- (disj);
  \draw[-Latex, thick] (fails) -- (tools);
  \draw[-Latex, thick] (disj) -- (tools);
\end{tikzpicture}
\caption{Why second-order MPEC theory departs from standard NLP theory.}
\end{figure}

\section{AVI-constrained mathematical programs}

We begin with the cleanest model. Consider
\begin{equation}\label{eq:avi}
\begin{aligned}
\min_{x,y}\quad & f(x,y)\\
\text{s.t.}\quad & (x,y)\in Z,\\
& Dx+Ey+b\le 0,\\
& (y'-y)^\T(Px+Qy+q)\ge 0
\quad
\forall y' \text{ such that } Dx+Ey'+b\le 0,
\end{aligned}
\end{equation}
where $Z$ is polyhedral.

This is an upper-level optimization problem constrained by an \emph{affine}
variational inequality. Under standard assumptions, the graph of the solution map
is piecewise polyhedral.

\subsection{Critical cone}

Let $\bar z=(\bar x,\bar y)$ be feasible and denote the feasible set by $\F$.
The MPEC critical cone is
\[
\C(\bar z;\F)
=
\left\{
dz\in T(\bar z;\F):
\nabla f(\bar z)^\T dz=0
\right\}.
\]

The essential fact is that because the lower-level set is polyhedral,
every critical direction can be realized by an actual feasible first-order arc
along one active polyhedral branch.

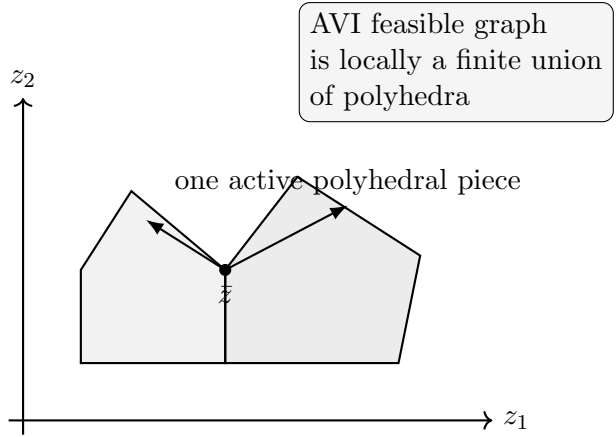
\begin{figure}[ht]
\centering
\begin{tikzpicture}[scale=0.95]
  \draw[axis] (-0.2,0) -- (6.5,0) node[right] {$z_1$};
  \draw[axis] (0,-0.2) -- (0,4.5) node[above] {$z_2$};

  \fill[gray!10] (0.8,0.8) -- (2.8,0.8) -- (2.8,2.1) -- (1.5,3.2) -- (0.8,2.1) -- cycle;
  \fill[gray!15] (2.8,0.8) -- (5.2,0.8) -- (5.5,2.3) -- (3.8,3.4) -- (2.8,2.1) -- cycle;

  \draw[thick] (0.8,0.8) -- (2.8,0.8) -- (2.8,2.1) -- (1.5,3.2) -- (0.8,2.1) -- cycle;
  \draw[thick] (2.8,0.8) -- (5.2,0.8) -- (5.5,2.3) -- (3.8,3.4) -- (2.8,2.1) -- cycle;

  \node[point, label=below:{$\bar z$}] at (2.8,2.1) {};
  \draw[ray] (2.8,2.1) -- (4.5,3.0) node[above] {one active polyhedral piece};
  \draw[ray] (2.8,2.1) -- (1.7,2.8);

  \node[labelbox] at (6.0,5.0) {
    AVI feasible graph\\
    is locally a finite union\\
    of polyhedra
  };
\end{tikzpicture}
\caption{Local union-of-polyhedra geometry in the AVI case.}
\end{figure}

\subsection{Clean second-order theory}

For AVI-constrained MPECs, the curvature of the lower-level feasible system is absent,
so the second-order condition reduces to the objective Hessian alone.

\begin{theorem}[Necessary condition in the polyhedral AVI case]
Let $\bar z$ be a local minimizer of \eqref{eq:avi} and assume $f$ is $C^2$
near $\bar z$. If the upper-level feasible region and the lower-level AVI data
are polyhedral, then
\[
dz^\T \nabla^2 f(\bar z)\, dz \ge 0
\qquad
\forall dz\in \C(\bar z;\F).
\]
That is, $\nabla^2 f(\bar z)$ is copositive on the MPEC critical cone.
\end{theorem}

\begin{theorem}[Sufficient condition in the polyhedral AVI case]
Assume $\bar z$ is stationary for \eqref{eq:avi}.
If
\[
dz^\T \nabla^2 f(\bar z)\, dz >0
\qquad
\forall dz\in \C(\bar z;\F)\setminus\{0\},
\]
then $\bar z$ is a strict local minimizer.
\end{theorem}

\begin{remark}
This is one of the most elegant features of the AVI case: there is no
constraint-curvature contribution in the second-order test.
The geometry is essentially piecewise linear.
\end{remark}

\section{NCP-constrained mathematical programs}

Now consider
\begin{equation}\label{eq:ncp-mpec}
\begin{aligned}
\min_{x,y}\quad & f(x,y)\\
\text{s.t.}\quad & x\in X,\\
& F(x,y)\ge 0,\qquad y\ge 0,\qquad y^\T F(x,y)=0,
\end{aligned}
\end{equation}
where $X\subset\R^n$ is polyhedral and $F:\R^{n+m}\to\R^m$ is $C^2$.

\subsection{Index-set decomposition}

At a feasible point $\bar z=(\bar x,\bar y)$ define
\begin{align*}
\alpha(\bar z) &= \{i : F_i(\bar z)=0 < \bar y_i\},\\
\beta(\bar z)  &= \{i : F_i(\bar z)=0 = \bar y_i\},\\
\gamma(\bar z) &= \{i : F_i(\bar z)>0 = \bar y_i\}.
\end{align*}

These sets encode the complementarity structure:

\begin{enumerate}[label=(\roman*)]
\item $\alpha$: active on the equation side $F_i=0$ with positive primal variable;
\item $\gamma$: active on the variable side $y_i=0$ with positive residual;
\item $\beta$: \emph{degenerate} indices where both are zero.
\end{enumerate}

The set $\beta(\bar z)$ is where all technical complications live.





\subsection{Critical cone}

The tangent directions $(dx,dy)$ satisfy
\[
dx\in T(\bar x;X),\qquad
dy_i=0\;\; (i\in\gamma(\bar z)),
\]
together with mixed complementarity linearization on $\beta(\bar z)$.
The critical cone is
\[
\C(\bar z;\F)
=
\{(dx,dy)\in T(\bar z;\F): \nabla f(\bar z)^\T(dx,dy)=0\}.
\]

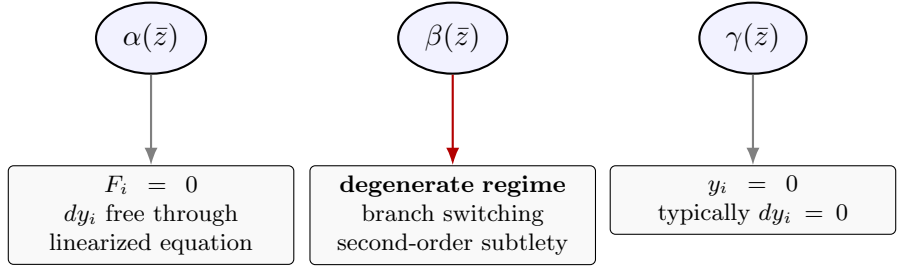
\begin{figure}[ht]
\centering
\begin{tikzpicture}[
    node distance=1.2cm and 2.5cm, 
    mybox/.style={draw, fill=gray!5, font=\footnotesize, align=center, rounded corners=2pt, text width=3.5cm},
    myellipse/.style={draw, ellipse, thick, align=center, inner sep=4pt, fill=blue!5}
  ]

  \node[myellipse] (alpha) {$\alpha(\bar z)$};
  \node[myellipse, right=of alpha] (beta) {$\beta(\bar z)$};
  \node[myellipse, right=of beta] (gamma) {$\gamma(\bar z)$};

  \node[mybox, below=of alpha] (a1) {
    $F_i=0$ \\ 
    $dy_i$ free through \\ 
    linearized equation
  };
  
  \node[mybox, below=of beta] (b1) {
    \textbf{degenerate regime} \\ 
    branch switching \\ 
    second-order subtlety
  };
  
  \node[mybox, below=of gamma] (g1) {
    $y_i=0$ \\ 
    typically $dy_i=0$
  };

  \draw[-Latex, thick, gray] (alpha) -- (a1);
  \draw[-Latex, thick, red!70!black] (beta) -- (b1); 
  \draw[-Latex, thick, gray] (gamma) -- (g1);

\end{tikzpicture}
\caption{Why $\beta(\bar z)$ is the technically difficult set in NCP-constrained MPECs.}
\end{figure}

\subsection{Multiplier-based second-order condition}

Because the lower-level data are nonlinear, curvature terms of the NCP mapping appear.
Formally, one expects a quadratic form of type
\[
dz^\T
\left(
\nabla^2 f(\bar z)-\sum_{i}\pi_i \nabla^2 F_i(\bar z)
\right)
dz,
\qquad dz=(dx,dy).
\]

\begin{remark}
This resembles the NLP Lagrangian Hessian, but the admissible multipliers
and admissible directions are more subtle. In particular, the sharp multiplier
set depends on the complementarity pattern inside $\beta(\bar z)$.
\end{remark}

A useful lecture-level summary is:

\begin{theorem}[Multiplier-form necessary condition, schematic]
Assume $\bar z$ is a local minimizer of \eqref{eq:ncp-mpec},
$F$ and $f$ are $C^2$, and the lower-level NCP is strongly regular at $\bar z$.
Then for every critical direction $dz=(dx,dy)$ satisfying the appropriate
strict-complementarity condition on degenerate indices, there exists an
MPEC multiplier $\pi$ such that
\[
dz^\T
\left(
\nabla^2 f(\bar z)-\sum_{i\in\alpha(\bar z)\cup\beta(\bar z)}
\pi_i \nabla^2 F_i(\bar z)
\right)
dz
\ge 0.
\]
\end{theorem}

\begin{theorem}[Multiplier-form sufficient condition, schematic]
Under the same regularity assumptions, if for every nonzero critical direction $dz$
and every admissible multiplier $\pi$ one has
\[
dz^\T
\left(
\nabla^2 f(\bar z)-\sum_{i\in\alpha(\bar z)\cup\beta(\bar z)}
\pi_i \nabla^2 F_i(\bar z)
\right)
dz
>0,
\]
then $\bar z$ is a strict local minimizer.
\end{theorem}

\section{Strong regularity and the implicit-programming viewpoint}

The most refined second-order theory comes from reducing the MPEC
to an optimization problem with a nonsmooth implicit lower-level response map.

\subsection{Strong regularity}

Fix $\bar z=(\bar x,\bar y)$ feasible for \eqref{eq:ncp-mpec}. Suppose the lower-level NCP
\[
F(x,y)\ge 0,\qquad y\ge 0,\qquad y^\T F(x,y)=0
\]
is \emph{strongly regular} at $(\bar x,\bar y)$.
Informally, strong regularity means that near $\bar x$, the solution $y$
is locally unique and depends on $x$ in a stable, piecewise-$C^1$ manner.

\begin{remark}
Strong regularity is the equilibrium analogue of a robust implicit-function property.
It is closely related to the local invertibility of a suitable generalized derivative
and is a central concept in Robinson's stability theory.
\end{remark}

Thus, near $\bar x$, one may write
\[
y = y(x),
\]
with $y(\cdot)$ a locally single-valued piecewise smooth map.

\begin{figure}[ht]
\centering
\begin{tikzpicture}[scale=1.2, >=Stealth]
  \tikzset{
    axis/.style={->, gray, thin},
    point/.style={circle, fill=black, inner sep=1.5pt},
    ray/.style={->, blue!70!black, thick},
    critical/.style={red!70!black, thick, dashed},
    labelbox/.style={draw, fill=white, font=\footnotesize, align=center, rounded corners=2pt}
  }

  \draw[axis] (-0.5,0) -- (6.5,0) node[right] {$x$};
  \draw[axis] (0,-0.5) -- (0,5.0) node[above] {$y(x)$};

  \coordinate (P) at (3.5,2.0);
  \draw[dashed, gray] (3.5,0) -- (P);
  \node[below] at (3.5,0) {$\bar{x}$};
  \node[point] at (P) {};
  \node[above left] at (P) {$(\bar{x}, \bar{y})$};

  \draw[thick] (1.0,1.2) .. controls (2.0,1.4) and (3.0,1.6) .. (P);
  \draw[thick] (P) .. controls (4.0,2.5) and (5.0,3.8) .. (5.8,4.5);
  
  \draw[ray] (P) -- (5.0,3.0) node[right, font=\small] {$y'(\bar{x}; dx)$};

  \node[red!70!black, font=\small] at (5.2,4.0) {piecewise-$C^1$};
  \draw[->, red!70!black, bend left=20] (4.8,3.8) to (4.3,2.8);

  \node[labelbox] at (2.0,4.2) {
    \textbf{Strong Regularity}\\
    $\Rightarrow$ local single-valued\\
    implicit solution map
  };

\end{tikzpicture}
\caption{Schematic local graph of the implicit solution map $y(x)$.}
\end{figure}
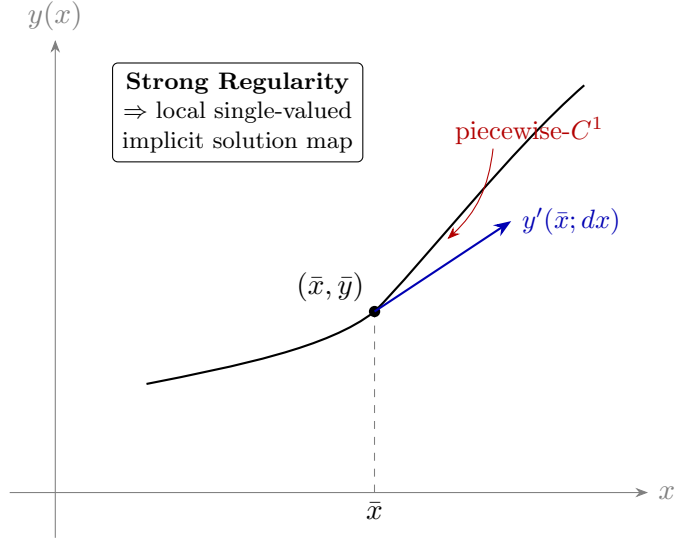

\subsection{Reduced formulation}

The MPEC becomes locally
\begin{equation}\label{eq:reduced}
\min_{x\in X} \varphi(x):=f(x,y(x)).
\end{equation}

This is no longer a smooth NLP in general, but it is sufficiently structured
to admit directional derivatives and, under strong regularity, a second-order directional expansion.

\subsection{First and second directional derivatives of the solution map}

Given a direction $dx$, the first directional derivative $y'(\bar x;dx)$
is characterized by a mixed linear complementarity problem.

Even more importantly, the second-order directional term
\[
y^{(2)}(\bar x;dx)
=
\lim_{\tau\downarrow 0}
\frac{y(\bar x+\tau dx)-\bar y-\tau y'(\bar x;dx)}{\tfrac12\tau^2}
\]
exists under the regularity assumptions and is characterized by a second mixed LCP.

\begin{figure}[ht]
\centering
\begin{tikzpicture}[scale=1.2, >=Stealth] 
  \tikzset{
    axis/.style={->, gray, thin},
    point/.style={circle, fill=black, inner sep=1.5pt},
    labelbox/.style={draw, fill=white, font=\footnotesize, align=center, rounded corners=2pt}
  }

  \draw[axis] (-0.5,0) -- (7.5,0) node[right] {$x$};
  \draw[axis] (0,-0.5) -- (0,5.5) node[above] {$y$};

  \draw[thick] (0.5,0.5) .. controls (2,1) and (3,1.5) .. (4,2.5) 
               node[pos=0.5, above left] {$y(x)$}
               .. controls (5,3.5) and (6,5) .. (6.5,5.5);

  \node[point, label=below right:{$(\bar{x}, \bar{y})$}] (p1) at (4,2.5) {};

  \coordinate (xtau) at (6.2,0);
  \draw[dashed] (6.2,0) -- (6.2,5.15); 
  \node[below] at (xtau) {$\bar{x} + \tau dx$};

  \draw[blue!70!black, thick] (2.5,1.3) -- (6.8,4.6) 
        node[pos=0.9, below right, font=\small] {First-order prediction};

  \node[point, fill=red!70!black] (p_actual) at (6.2,5.15) {};
  \node[right=3pt, red!70!black, font=\small] at (p_actual) {$y(\bar{x} + \tau dx)$};

  \draw[<->, thick, red!70!black] (6.2,4.15) -- (6.2,5.15) 
        node[midway, left, font=\footnotesize] {$\frac{1}{2}\tau^2 y^{(2)}(\bar{x}; dx)$};

  \node[labelbox] at (2.5,4.5) {
    \textbf{Second-order term}\\
    Measures the quadratic\\
    deviation from the\\
    linear approximation
  };

\end{tikzpicture}
\caption{Interpretation of $y^{(2)}(\bar x;dx)$ as second-order deviation.}
\end{figure}
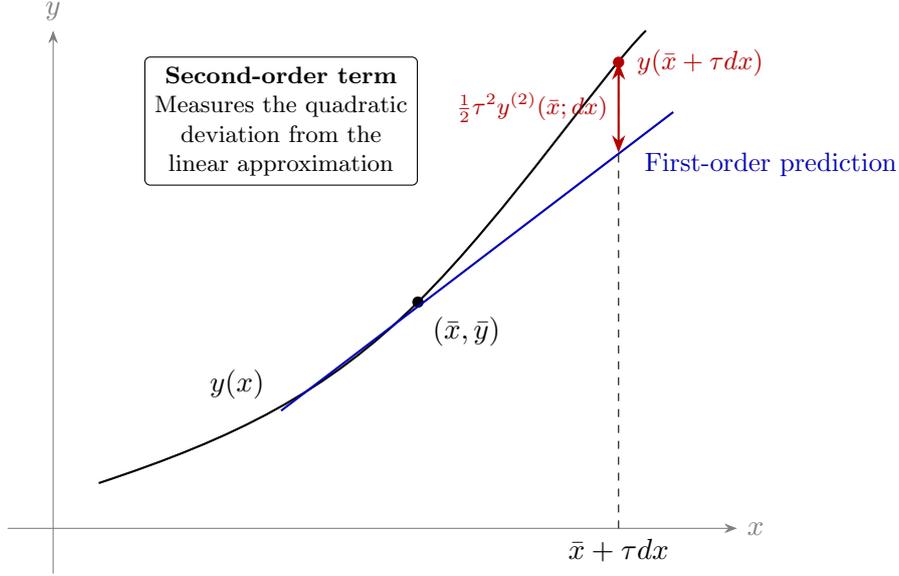

\begin{remark}
This object is the correct replacement for the ordinary Hessian of an implicit $C^2$ map.
When the degenerate index set $\beta(\bar z)$ is empty, it reduces to the usual second derivative
computed from the classical implicit-function theorem.
\end{remark}

\subsection{Second-order conditions in reduced form}

The reduced second-order expression is
\[
dz^\T \nabla^2 f(\bar z)\, dz
+
\nabla_y f(\bar z)^\T y^{(2)}(\bar x;dx),
\qquad
dz=(dx,y'(\bar x;dx)).
\]

\begin{theorem}[Implicit-programming necessary condition]
Assume $\bar z=(\bar x,\bar y)$ is a stationary point of \eqref{eq:ncp-mpec},
the lower-level NCP is strongly regular at $\bar z$,
and $f,F$ are $C^2$ near $\bar z$.
If $\bar z$ is a local minimizer, then for every critical direction
$dz=(dx,dy)\in \C(\bar z;\F)$ with $dy=y'(\bar x;dx)$,
\[
dz^\T \nabla^2 f(\bar z)\, dz
+
\nabla_y f(\bar z)^\T y^{(2)}(\bar x;dx)
\ge 0.
\]
\end{theorem}

\begin{theorem}[Implicit-programming sufficient condition]
Under the same assumptions, if
\[
dz^\T \nabla^2 f(\bar z)\, dz
+
\nabla_y f(\bar z)^\T y^{(2)}(\bar x;dx)
>0
\]
for every nonzero critical direction $dz=(dx,y'(\bar x;dx))$, then $\bar z$
is a strict local minimizer.
\end{theorem}

\section{From implicit-programming back to multipliers}

The reduced condition can be converted back into a multiplier statement.
The conversion is nontrivial because the second-order directional term
$y^{(2)}(\bar x;dx)$ itself satisfies a mixed LCP whose multipliers depend on the active degenerate pattern.

The outcome is a sharpened multiplier formula:
\[
dz^\T
\left(
\nabla^2 f(\bar z)-\sum_i \pi_i \nabla^2 F_i(\bar z)
\right)dz
+
\text{additional nonnegative terms}
\ge 0.
\]

The additional nonnegative terms quantify second-order contributions from
degenerate indices in $\beta(\bar z)$. They explain precisely why the naive
multiplier form from classical NLP is incomplete in the MPEC setting.

\begin{figure}[ht]
\centering
\begin{tikzpicture}[node distance=4.5cm]
  \node[mybox] (lhs) {Reduced expression\\
  $dz^\T\nabla^2 f(\bar z)dz+\nabla_y f(\bar z)^\T y^{(2)}(\bar x;dx)$};
  \node[mybox, right=of lhs] (rhs) {Multiplier form\\
  $dz^\T\!\left(\nabla^2 f-\sum_i \pi_i\nabla^2 F_i\right)\!dz$\\
  $+\;$nonnegative switching terms};

  \draw[-Latex, thick] (lhs) -- node[above] {eliminate $y^{(2)}$ via LCP} (rhs);
\end{tikzpicture}
\caption{The refined multiplier formula is obtained by eliminating $y^{(2)}$ through its characterization system.}
\end{figure}
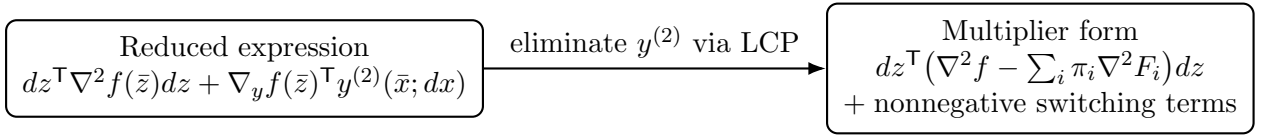

\section{KKT-constrained MPECs}

Consider now an MPEC whose lower-level equilibrium is expressed through KKT conditions:
\begin{equation}\label{eq:kkt-mpec}
\begin{aligned}
\min_{x,y,\lambda}\quad & f(x,y)\\
\text{s.t.}\quad & Gx+a\le 0,\\
& 0 = F(x,y)+\nabla_y g(x,y)^\T \lambda,\\
& g(x,y)\le 0,\qquad \lambda\ge 0,\qquad \lambda^\T g(x,y)=0.
\end{aligned}
\end{equation}

Here $\lambda$ is the lower-level multiplier.

\subsection{MPEC Lagrangian}

The natural MPEC Lagrangian is
\[
\Lag_{\MPEC}(x,y,\lambda;\zeta,\pi,\eta)
=
f(x,y)
+\zeta^\T(Gx+a)
-\pi^\T\big(F(x,y)+\nabla_y g(x,y)^\T\lambda\big)
+\eta^\T g(x,y).
\]

The second-order tests are now written in the enlarged variable $w=(x,y,\lambda)$.

\begin{figure}[ht]
\centering
\begin{tikzpicture}[node distance=2.1cm]
  \node[mybox] (x) {$x$\\upper variable};
  \node[mybox, right=of x] (y) {$y$\\primal lower variable};
  \node[mybox, right=of y] (lam) {$\lambda$\\lower multiplier};

  \node[mybox, below=of y] (kkt) {lower-level KKT system\\
  $0=F(x,y)+\nabla_y g(x,y)^\T\lambda$\\
  $g(x,y)\le 0,\;\lambda\ge 0,\;\lambda^\T g(x,y)=0$};

  \node[mybox, below=of lam] (lag) {$\Lag_{\MPEC}(x,y,\lambda;\zeta,\pi,\eta)$};

  \draw[-Latex, thick] (x) -- (kkt);
  \draw[-Latex, thick] (y) -- (kkt);
  \draw[-Latex, thick] (lam) -- (kkt);
  \draw[-Latex, thick] (kkt) -- (lag);
\end{tikzpicture}
\caption{Variable roles in the KKT-constrained MPEC formulation.}
\end{figure}
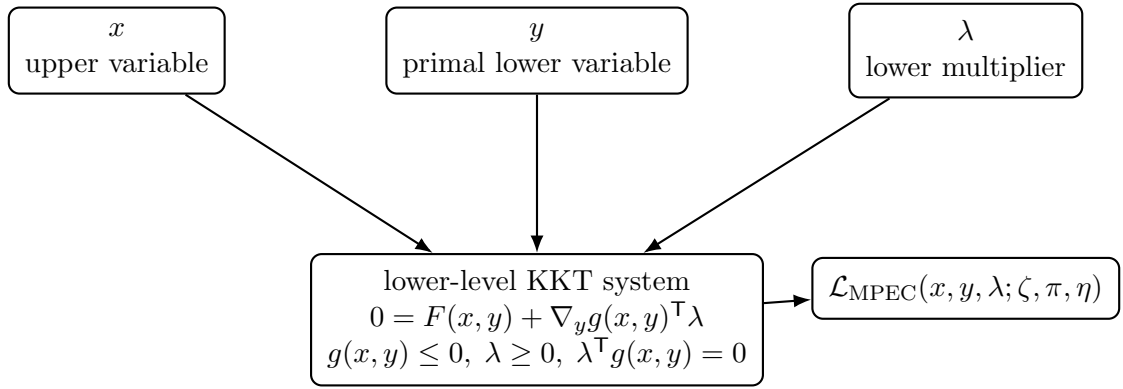

\subsection{Structure of the sufficient condition}

Under LICQ-type conditions and strong regularity of the lower-level KKT system,
one obtains a result of the form:

\begin{theorem}[Schematic KKT-based sufficient condition]
Assume $(\bar x,\bar y,\bar\lambda)$ is stationary for \eqref{eq:kkt-mpec},
the lower-level KKT system is strongly regular,
and for every nonzero critical direction $dw$ there exist admissible multipliers
$(\zeta,\pi,\eta)$ such that
\[
dw^\T \nabla^2_{ww}\Lag_{\MPEC}(\bar w;\zeta,\pi,\eta)\, dw >0.
\]
Then $\bar w$ is a strict local minimizer of the KKT-constrained MPEC.
\end{theorem}

\begin{remark}
A notable feature is that strict local minimality in the enlarged space often implies
local uniqueness of the lower-level KKT multiplier. This is consistent with
the general philosophy that strong second-order growth suppresses multiplier nonuniqueness.
\end{remark}

\section{Piecewise-programming approach}

A third route is to split complementarity explicitly into smooth regimes.

For the NCP model \eqref{eq:ncp-mpec}, each degenerate index $i\in\beta(\bar z)$
admits two local branches:
\[
F_i(x,y)=0
\quad \text{or} \quad
y_i=0.
\]
Choosing one branch for each degenerate index yields a smooth NLP subproblem.
There are finitely many such subproblems.

\subsection{Smooth pieces}

Fix a subset $\mathcal{B}\subseteq \beta(\bar z)$.
Interpret indices in $\mathcal{B}$ via the branch $F_i=0$ and those in $\beta(\bar z)\setminus \mathcal{B}$
via the branch $y_i=0$. This produces a smooth NLP on one active piece.

On each piece, standard NLP second-order theory applies directly.

\begin{figure}[ht]
\centering
\begin{tikzpicture}[scale=1.0]
  \draw[axis] (-0.2,0) -- (6.7,0) node[right] {$z_1$};
  \draw[axis] (0,-0.2) -- (0,4.6) node[above] {$z_2$};

  \fill[gray!10] (0.8,0.8) -- (3.2,0.8) -- (3.2,2.6) -- cycle;
  \draw[thick] (0.8,0.8) -- (3.2,0.8) -- (3.2,2.6) -- cycle;
  \node at (2.45,1.45) {$\mathcal P_1$};

  \fill[gray!18] (3.2,0.8) -- (5.8,0.8) -- (3.2,2.6) -- cycle;
  \draw[thick] (3.2,0.8) -- (5.8,0.8) -- (3.2,2.6) -- cycle;
  \node at (4.3,1.35) {$\mathcal P_2$};

  \node[point, label=above left:{$\bar z$}] at (3.2,2.6) {};

  \node[labelbox] at (5.1,3.8) {
    complementarity degeneracy\\
    is decomposed into\\
    finitely many NLP branches
  };
\end{tikzpicture}
\caption{Piecewise-programming viewpoint: solve second-order analysis branch by branch.}
\end{figure}
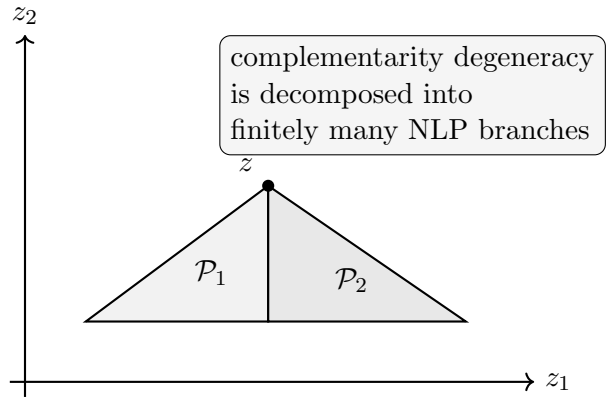

\subsection{Advantages and drawbacks}

\begin{enumerate}[label=(\roman*)]
\item Advantage: conceptually elementary; one can literally apply the classical KKT and second-order theory to each piece.
\item Drawback: potentially combinatorial; one must inspect many branch patterns.
\item Interpretation: the piecewise approach explains why complementarity constraints act
``as if linear'' inside each branch, even though the global feasible set is nonsmooth.
\end{enumerate}

\begin{theorem}[Schematic piecewise sufficient condition]
Suppose that for every relevant complementarity branch, the corresponding smooth NLP
satisfies MFCQ and its Lagrangian Hessian is strictly positive on the associated critical cone.
Then the original MPEC has a strict local minimizer at the reference point.
\end{theorem}

\section{Worked micro-example}

Consider
\begin{equation}\label{eq:micro}
\begin{aligned}
\min_{x,y}\quad & \frac12(x^2-y^2)\\
\text{s.t.}\quad & x\ge 0,\; y\ge 0,\; y(y^4+y+x)=0,\; y^4+y+x\ge 0.
\end{aligned}
\end{equation}

At $(\bar x,\bar y)=(0,0)$ the feasible set is locally
\[
\F = \R_+\times\{0\},
\]
so the tangent cone is
\[
T(\bar z;\F)=\R_+\times\{0\}.
\]
The gradient is zero at the origin, hence
\[
\C(\bar z;\F)=T(\bar z;\F).
\]
The Hessian of the objective is
\[
\nabla^2 f(0,0)=
\begin{pmatrix}
1 & 0\\
0 & -1
\end{pmatrix}.
\]
Therefore, for any critical direction $dz=(d_x,0)$,
\[
dz^\T \nabla^2 f(0,0)\, dz = d_x^2\ge 0,
\]
with strict inequality for nonzero $d_x$.
Hence the second-order growth condition holds on the critical cone,
and $(0,0)$ is a strict local minimizer.





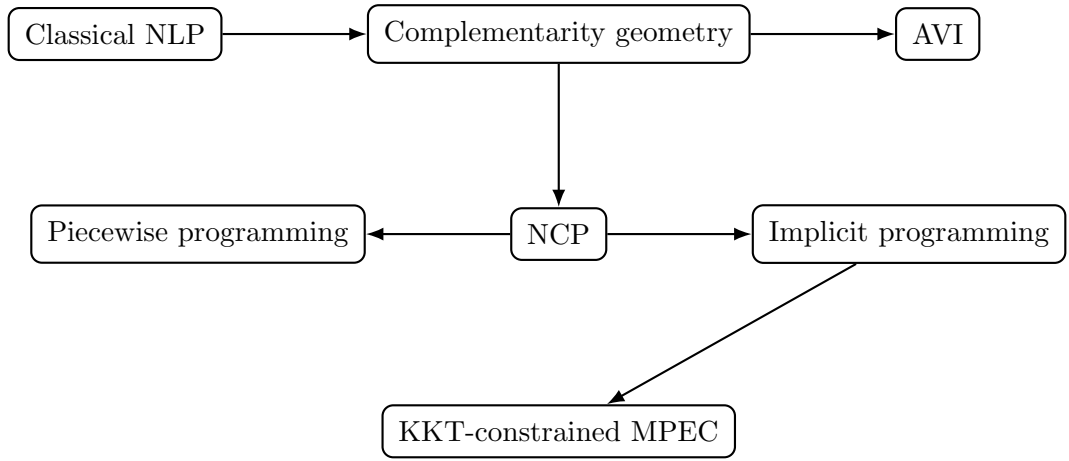
\begin{figure}[ht]
\centering
\begin{tikzpicture}[node distance=1.9cm]
  \node[mybox] (a) {Classical NLP};
  \node[mybox, right=of a] (b) {Complementarity geometry};
  \node[mybox, right=of b] (c) {AVI};
  \node[mybox, below=of b] (d) {NCP};
  \node[mybox, left=of d] (e) {Piecewise programming};
  \node[mybox, right=of d] (f) {Implicit programming};
  \node[mybox, below=of d] (g) {KKT-constrained MPEC};

  \draw[-Latex, thick] (a) -- (b);
  \draw[-Latex, thick] (b) -- (c);
  \draw[-Latex, thick] (b) -- (d);
  \draw[-Latex, thick] (d) -- (e);
  \draw[-Latex, thick] (d) -- (f);
  \draw[-Latex, thick] (f) -- (g);
\end{tikzpicture}
\caption{A good pedagogical sequencing for reading the notes.}
\end{figure}

From a research perspective, the main moral is:

\begin{quote}
Second-order optimality for MPECs is not merely ``NLP plus complementarity.''
It is a synthesis of variational geometry, generalized differentiation,
stability theory, and branchwise smooth analysis.
\end{quote}

\section{Exercises}
\subsection*{True or False Questions}

\begin{enumerate}[label=(\roman*)]
\item Standard nonlinear programming (NLP) second-order optimality theory, such as Theorem 5.1.1, can be directly applied to an Affine Variational Inequality (AVI) constrained Mathematical Program because standard constraint qualifications like the Mangasarian-Fromovitz Constraint Qualification (MFCQ) inherently hold for MPECs.

\item When formulating the second-order optimality conditions for an AVI constrained Mathematical Program, the second derivatives of the nonlinear complementarity constraint $\lambda^T(Dx+Ey+b)=0$ are omitted from the Hessian matrix of the upper-level objective function.

\item For a Mathematical Program with Affine Equilibrium Constraints (MPAEC) featuring a quadratic upper-level objective function, the copositivity of the objective function's Hessian matrix on the critical cone is a necessary and sufficient condition for a stationary point to be a local minimum.

\item  In the context of an NCP constrained Mathematical Program, the second-order necessary optimality condition (Theorem 5.3.1) guarantees that the restricted copositivity condition holds for all vectors within the critical cone $\mathcal{C}(\bar{z};\mathcal{F})$.

\item The second-order directional derivative of the implicit solution function, $y^{(2)}(\bar{x}; dx)$, is positively homogeneous of degree 2, meaning that $y^{(2)}(\bar{x}; \tau dx) = \tau^2 y^{(2)}(\bar{x}; dx)$ for all $\tau \ge 0$.

\item By utilizing the implicit programming approach and evaluating the second-order directional derivative $y^{(2)}(\bar{x}; dx)$, it is possible to formulate a second-order necessary condition for the NCP constrained MP that completely avoids the use of MPEC multipliers from the set $P(\bar{z})$.

\item In a KKT Constrained Mathematical Program, if the second-order sufficiency condition holds (i.e., the MPEC Lagrangean's Hessian is strictly copositive on the critical cone), it implies that the KKT multiplier $\bar{\lambda}$ associated with the solution is strictly non-unique.

\item If the inner Nonlinear Complementarity Problem (NCP) is defined by purely affine functions, the second-order directional derivative $y^{(2)}(\bar{x}; dx)$ evaluates to zero for all critical vectors $dz \in \mathcal{C}(\bar{z};\mathcal{F})$.

\item In the piecewise programming approach, establishing the Linear Independence Constraint Qualification (LICQ) for a disjunctive subproblem guarantees that the Strong Regularity Condition (SRC) holds for the overall NCP constrained Mathematical Program.

\item The critical cone of the relaxed NLP, $\mathcal{C}_{rl}(\bar{w})$, is always a subset of the union of the critical cones $\mathcal{C}_{\alpha}(\bar{w})$ derived from the restricted nonlinear programs corresponding to various index subsets $\alpha$.
\end{enumerate}

\subsection*{Computational Questions}

\subsubsection*{Question 1. Explicit Computation of the MPEC Lagrangean and SCOC}
    Consider the following NCP constrained Mathematical Program in 2 variables, common in power-control equilibrium models:
    $$\text{minimize} \quad f(x,y) = \frac{1}{2}(x^2 - y^2)$$
    $$\text{subject to} \quad x \ge 0, \ y \ge 0$$
    $$F(x,y) = y^4 + y + x \ge 0$$
    $$y F(x,y) = 0$$
    Let $\overline{z} \equiv (\overline{x},\overline{y}) = (0,0)$.
    \begin{enumerate}[label=(\roman*)]
        \item Compute the gradients $\nabla f(\overline{z})$ and $\nabla F(\overline{z})$, and their respective Hessian matrices $\nabla^2 f(\overline{z})$ and $\nabla^2 F(\overline{z})$.
        \item Determine the index sets $\alpha(\overline{z})$, $\beta(\overline{z})$, and $\gamma(\overline{z})$. 
        \item By hand, compute the tangent cone $\mathcal{T}(\overline{z};\mathcal{F})$ and the critical cone $\mathcal{C}(\overline{z};\mathcal{F})$. 
        \item Using Theorem 5.3.1, verify the second-order necessary optimality condition at $\overline{z}$ by explicitly writing out the restricted copositivity condition:
        $$dz^T(\nabla^2 f(\overline{z}) - \pi \nabla^2 F(\overline{z}))dz \ge 0$$
        Demonstrate algebraically whether this holds for all critical vectors $dz \in \mathcal{C}(\overline{z};\mathcal{F})$ and valid multipliers $\pi$.
    \end{enumerate}

\subsubsection*{Question 2. Strong Regularity Condition (SRC) Computation}
    In network optimization, ensuring the stability of the inner NCP is critical. Consider an inner NCP parameterized by $x \in \mathbb{R}^2$, defined by the function:
    $$F(x,y) = \begin{pmatrix} x_1^2 + y_1 \\ x_2^2 + y_2 \end{pmatrix}$$
    Let $\overline{z} \equiv (\overline{x}_1, \overline{x}_2, \overline{y}_1, \overline{y}_2) = (0,0,0,0)$. 
    \begin{enumerate}[label=(\roman*)]
        \item Identify the active index set $\beta(\overline{z})$.
        \item Compute the partial Jacobian matrices $\nabla_x F(\overline{z})$ and $\nabla_y F(\overline{z})$.
        \item The Strong Regularity Condition (SRC) at $\overline{z}$ for an NCP requires that the partial Jacobian $\nabla_{\alpha(\overline{z})}F_{\alpha(\overline{z})}(\overline{z})$ is nonsingular and a specific Schur complement is a P-matrix. Since $\alpha(\overline{z}) = \emptyset$ here, the condition reduces strictly to the properties of $\nabla_{\beta(\overline{z})}F_{\beta(\overline{z})}(\overline{z})$. Compute this matrix and prove whether the SRC holds at $\overline{z}$.
    \end{enumerate}

\subsubsection*{Question 3. Implicit Programming and Second-Order Directional Derivatives}
    Consider a scalar parameter $x \in \mathbb{R}$ and a primary variable $y \in \mathbb{R}$ defining a parametric NCP:
    $$F(x,y) = x + y \ge 0, \quad y \ge 0, \quad y(x+y) = 0$$
    Assume we are at the stationary point $\overline{x} = 0, \overline{y} = 0$. 
    \begin{enumerate}[label=(\roman*)]
        \item For a given perturbation direction $dx = 1$, compute the first-order directional derivative $y'(\overline{x}; dx)$ by solving the corresponding mixed LCP.
        \item Define the index sets $d\beta_a(\overline{x};dx)$, $d\beta_b(\overline{x};dx)$, and $d\beta_g(\overline{x};dx)$ for $dx = 1$. Determine which set contains the index $i=1$.
        \item Formulate the mixed LCP system for the second-order directional derivative $y^{(2)}(\overline{x}; dx)$ as given in Proposition 5.4.1(b). Note that the system depends on $dz^T \nabla^2 F_i(\overline{z}) dz$. 
        \item Solve the mixed LCP by hand to find the exact numerical value of $y^{(2)}(\overline{x}; dx)$ for $dx = 1$.
    \end{enumerate}

\subsection*{Proof Questions}

\subsubsection*{Question 1. AVI Constrained Mathematical Programs}
    Consider a Mathematical Program with Affine Equilibrium Constraints (MPAEC) where the upper-level objective function $f(x,y)$ is twice continuously differentiable and the constraint set $Z$ is polyhedral. 
    \begin{enumerate}[label=(\roman*)]
        \item Explain geometrically and analytically why the standard Mangasarian-Fromovitz Constraint Qualification (MFCQ) generally fails for the feasible region $\mathcal{F}$ of an MPAEC.
        \item When formulating the second-order necessary optimality conditions for the MPAEC, the second derivatives of the nonlinear complementarity constraint $\lambda^T(Dx+Ey+b)=0$ are excluded from the MPEC Lagrangean's Hessian. Justify this exclusion using the piecewise linear structure of the tangent map $\mathcal{LS}_{(\overline{z},\lambda)}$.
    \end{enumerate}

\subsubsection*{Question 2. NCP Constraints and Strong Regularity}
    Consider an NCP constrained MP at a stationary point $\overline{z} \equiv (\overline{x},\overline{y}) \in \mathcal{F}$. Assume the Strong Regularity Condition (SRC) holds at $\overline{z}$.
    \begin{enumerate}[label=(\roman*)]
        \item Define the SRC in terms of the index sets $\alpha(\overline{z})$ and $\beta(\overline{z})$.
        \item Theorem 5.3.1 establishes a second-order necessary condition $\forall dz \in \mathcal{C}(\overline{z};\mathcal{F})$, but it restricts the critical vectors $dz \equiv (dx,dy)$ to those satisfying the strict complementarity condition:
        $$dy_i + \nabla F_i(\overline{z})^T dz > 0, \quad \forall i \in \beta(\overline{z})$$
        Derive why this specific restriction is necessary when utilizing the standard multiplier-based approach, and explain the limitations it imposes on the general second-order theory.
    \end{enumerate}

\subsubsection*{Question 3. Implicit Programming \& Directional Derivatives}
    To overcome the strict complementarity restriction in Question 2, we turn to the implicit programming formulation. Let $y(x)$ be the $PC^1$ implicit solution function of the parametric NCP. 
    \begin{enumerate}[label=(\roman*)]
        \item State the limit definition of the second-order directional derivative $y^{(2)}(\overline{x}; dx)$.
        \item Prove that $y^{(2)}(\overline{x}; \cdot)$ is positively homogeneous of degree 2; that is, show:
        $$y^{(2)}(\overline{x}; \tau dx) = \tau^2 y^{(2)}(\overline{x}; dx), \quad \forall \tau \ge 0$$
        \item Let $\overline{z}$ be a strongly regular local minimum. Prove the second-order necessary condition:
        $$dz^T \nabla^2 f(\overline{z}) dz + \nabla_y f(\overline{z})^T y^{(2)}(\overline{x}; dx) \ge 0$$
        for all critical vectors $dz \in \mathcal{C}(\overline{z};\mathcal{F})$. How does this formulation successfully avoid dependence on the MPEC multipliers $\pi \in P(\overline{z})$?
    \end{enumerate}

\subsubsection*{Question 4. Piecewise Programming and the Relaxed NLP}
    Consider a KKT constrained Mathematical Program and a strictly complementary local minimum $\overline{w} \equiv (\overline{z},\overline{\lambda}) \in \mathcal{F}^{KKT}$.
    \begin{enumerate}[label=(\roman*)]
        \item Formulate the relaxed nonlinear program (NLP) associated with $\overline{w}$. Define its critical cone, $\mathcal{C}_{rl}(\overline{w})$.
        \item Assume the MFCQ holds for the relaxed NLP at $\overline{w}$. Prove that if the MPEC Lagrangean's Hessian matrix, $\nabla_{ww}^2\mathcal{L}^{MPEC}(\overline{w},\zeta,\pi,\eta)$, is strictly copositive on the relaxed critical cone $\mathcal{C}_{rl}(\overline{w})$, then the KKT multiplier $\overline{\lambda}$ is uniquely determined (i.e., $M(\overline{z}) = \{\overline{\lambda}\}$).
    \end{enumerate}
  
\section{Solutions}
\subsection*{Solutions for True or False}
\begin{enumerate}[label=(\roman*)]
    \item     False. The MFCQ generally fails to hold for MPECs due to the special geometric features of the tangent cone. Therefore, standard NLP theorems like Theorem 5.1.1 cannot be directly applied without utilizing alternative approaches (such as piecewise or implicit programming).
    \item True. Because the constraints define a piecewise linear mixed complementarity problem, the problem can be viewed as an NLP with disjunctive linear constraints. Therefore, the second derivatives of the complementarity condition do not contribute to the Hessian of the Lagrangean.
    \item True. Theorem 5.2.4(a) extends classical quadratic programming facts, stating that for an MPAEC with a quadratic objective and polyhedral constraints, stationarity plus copositivity of the Hessian on the critical cone is both necessary and sufficient for a local minimum.
    \item False. A major limitation (or "imperfection," as the text notes) of Theorem 5.3.1 is that the condition does not hold for all vectors in the critical cone. It only applies to those critical vectors $dz$ that also satisfy a strict complementarity condition (specifically, $dy_i + \nabla F_i(\bar{z})^T dz > 0$ for all $i \in \beta(\bar{z})$).
    \item True. Proposition 5.4.1(c) explicitly proves this homogeneity property based on the fact that the Hessian terms are positively homogeneous of degree 2 in $dx$.
    \item True. Theorem 5.4.2(a) provides a second-order necessary condition strictly in terms of the Hessian of the objective function and the second-order directional derivative $y^{(2)}(\bar{x}; dx)$, independent of the multiplier set $P(\bar{z})$.
    \item False. It implies the exact opposite: the uniqueness of the KKT multiplier $\bar{\lambda}$. This is due to a special structure of the MPEC Lagrangean function, which contains no quadratic term in $\lambda$ (its second partial derivative with respect to $\lambda$ is identically zero).
    \item True. As discussed under Corollary 5.4.5, if the constraint functions $F_i(x,y)$ are affine, their second derivatives (Hessians) are zero. By the uniqueness of the solution to the mixed LCP system, $y^{(2)}(\bar{x}; dx)$ must evaluate to zero.
    \item False. Example 5.6.3 explicitly demonstrates that LICQ for a subproblem and SRC for the NCP constrained MP are not comparable. LICQ does not imply SRC (because multiple solutions can violate the B-Implicit Function condition), and SRC does not imply LICQ.
    \item False. The relationship is the reverse. As shown in equation (12) of the text, the constraints of the relaxed NLP are a relaxation of the restricted subproblems. Thus, the union of the restricted critical cones $\mathcal{C}_{\alpha}(\bar{w})$ is a subset of the relaxed critical cone $\mathcal{C}_{rl}(\bar{w})$.
\end{enumerate}

\subsection*{Solutions for computational questions}

\paragraph{Problem 1: Explicit NCP Constrained MP}

(a) Gradients and Hessians

First, we take the derivatives of the upper-level objective function $f(x,y) = \frac{1}{2}(x^2 - y^2)$:

$\nabla f(x,y) = \begin{pmatrix} x \\ -y \end{pmatrix} \implies \nabla f(\overline{z}) = \begin{pmatrix} 0 \\ 0 \end{pmatrix}$

$\nabla^2 f(x,y) = \begin{pmatrix} 1 & 0 \\ 0 & -1 \end{pmatrix} \implies \nabla^2 f(\overline{z}) = \begin{pmatrix} 1 & 0 \\ 0 & -1 \end{pmatrix}$

Next, for the constraint function $F(x,y) = y^4 + y + x$:

$\nabla F(x,y) = \begin{pmatrix} 1 \\ 4y^3 + 1 \end{pmatrix} \implies \nabla F(\overline{z}) = \begin{pmatrix} 1 \\ 1 \end{pmatrix}$

$\nabla^2 F(x,y) = \begin{pmatrix} 0 & 0 \\ 0 & 12y^2 \end{pmatrix} \implies \nabla^2 F(\overline{z}) = \begin{pmatrix} 0 & 0 \\ 0 & 0 \end{pmatrix}$

(b) We calculate index sets. At $\overline{z} = (0,0)$, we have $F(\overline{z}) = 0$ and $\overline{y} = 0$.
Based on the definitions in the text:

$\alpha(\overline{z}) = \{i : F_i(\overline{z}) = 0, \overline{y}_i > 0\} = \emptyset$

$\beta(\overline{z}) = \{i : F_i(\overline{z}) = 0, \overline{y}_i = 0\} = \{1\}$ (since there is only 1 constraint)

$\gamma(\overline{z}) = \{i : F_i(\overline{z}) > 0, \overline{y}_i = 0\} = \emptyset$

(c) We calculate tangent and critical cones. The feasible region $\mathcal{F}$ is defined by $x \ge 0$, $y \ge 0$, $y^4+y+x \ge 0$, and $y(y^4+y+x) = 0$.
Since $x \ge 0$ and $y \ge 0$, the term $(y^4+y+x)$ is strictly positive for any $y > 0$. Therefore, to satisfy the complementarity constraint $y(y^4+y+x) = 0$, we must strictly have $y = 0$.
This simplifies the feasible region to the non-negative x-axis: $\mathcal{F} = \{(x,0) : x \ge 0\}$.
The tangent cone at the origin is exactly the region itself: $\mathcal{T}(\overline{z};\mathcal{F}) = \{(dx, dy) : dx \ge 0, dy = 0\}$.
The critical cone is $\mathcal{C}(\overline{z};\mathcal{F}) = \mathcal{T}(\overline{z};\mathcal{F}) \cap \nabla f(\overline{z})^\perp$. Since $\nabla f(\overline{z}) = 0$, its orthogonal complement is all of $\mathbb{R}^2$.
Therefore, $\mathcal{C}(\overline{z};\mathcal{F}) = \{(dx,0) : dx \ge 0\}$.

(d) We verify Theorem 5.3.1 (SCOC).
For a critical vector $dz = (dx, 0)$ with $dx \ge 0$:
$$dz^T(\nabla^2 f(\overline{z}) - \pi \nabla^2 F(\overline{z}))dz \ge 0$$
Since $\nabla^2 F(\overline{z}) = 0$ (the zero matrix), the multiplier $\pi$ drops out completely. We only need to evaluate:
$$dz^T \nabla^2 f(\overline{z}) dz = \begin{pmatrix} dx & 0 \end{pmatrix} \begin{pmatrix} 1 & 0 \\ 0 & -1 \end{pmatrix} \begin{pmatrix} dx \\ 0 \end{pmatrix} = dx^2$$
Since $dx^2 \ge 0$ for all real $dx$, the condition $dx^2 \ge 0$ holds trivially for all critical vectors $dz \in \mathcal{C}(\overline{z};\mathcal{F})$.

\paragraph{Problem 2: Checking the Strong Regularity Condition (SRC)}

(a) We calculate the active index set $\beta(\overline{z})$.
At $\overline{z} = (0,0,0,0)$, we calculate $F(\overline{z})$:
$F_1(0,0) = 0^2 + 0 = 0$ and $F_2(0,0) = 0^2 + 0 = 0$.
Since $\overline{y}_1 = 0$ and $\overline{y}_2 = 0$, both components satisfy $F_i = 0$ and $y_i = 0$.
Thus, $\beta(\overline{z}) = \{1, 2\}$.

(b) We calculate partial Jacobian matrices.
Taking the derivatives of $F(x,y)$ with respect to $x$ and $y$:
$\nabla_x F(x,y) = \begin{pmatrix} 2x_1 & 0 \\ 0 & 2x_2 \end{pmatrix} \implies \nabla_x F(\overline{z}) = \begin{pmatrix} 0 & 0 \\ 0 & 0 \end{pmatrix}$
$\nabla_y F(x,y) = \begin{pmatrix} 1 & 0 \\ 0 & 1 \end{pmatrix} \implies \nabla_y F(\overline{z}) = \begin{pmatrix} 1 & 0 \\ 0 & 1 \end{pmatrix} = I_2$

(c) We prove the SRC.
Because $\alpha(\overline{z}) = \emptyset$, the SRC definition simply requires that the matrix $\nabla_{\beta(\overline{z})} F_{\beta(\overline{z})}(\overline{z})$ is a P-matrix.
Here, $\nabla_{\beta(\overline{z})} F_{\beta(\overline{z})}(\overline{z}) = \nabla_y F(\overline{z}) = \begin{pmatrix} 1 & 0 \\ 0 & 1 \end{pmatrix}$.
A matrix is a P-matrix if and only if all of its principal minors are strictly positive.
The principal minors of the $2 \times 2$ identity matrix are as follows.
$1 \times 1$ minor is $\det([1]) = 1 > 0$.
$2 \times 2$ minor is $\det(I_2) = 1 > 0$.
Because all principal minors are strictly positive, it is a P-matrix. Therefore, the Strong Regularity Condition (SRC) holds at $\overline{z}$.

\paragraph{Problem 3: Implicit Programming and Second-Order Directional Derivatives}

(a) We calculate first-order directional derivative $y'(\overline{x}; dx)$.
The mixed LCP to find $dy = y'(\overline{x}; 1)$ is formulated using the linearization:
$$\nabla_x F(\overline{z}) dx + \nabla_y F(\overline{z}) dy = (1)(1) + (1)(dy) = 1 + dy.$$
The LCP for the single active index $i \in \beta$ are
$1 + dy \ge 0$,
$dy \ge 0$, and
$dy(1 + dy) = 0$.
If we assume $dy > 0$, then $1 + dy > 1$. Their product would be strictly positive, violating the complementarity condition (3). Therefore, we must have $dy = 0$.
So, $y'(\overline{x}; 1) = 0$.

(b) We calculate index sets.
We evaluate the terms for the index sets using $dx = 1$ and $d\overline{y} = 0$:
$\nabla_x F dx + \nabla_y F d\overline{y} = 1 + 0 = 1$, and
$d\overline{y} = 0$.
Checking the definitions gives

$d\beta_a$: Requires $\nabla F \cdot dz = 0$ and $d\overline{y} > 0$ (False).

$d\beta_b$: Requires $\nabla F \cdot dz = 0$ and $d\overline{y} = 0$ (False).

$d\beta_g$: Requires $\nabla F \cdot dz > 0$ and $d\overline{y} = 0$. We have $1 > 0$ and $0 = 0$ (True).

Therefore, index $i=1 \in d\beta_g(\overline{x}; 1)$.

(c) We formulate the mixed LCP for $y^{(2)}$.
As established in Proposition 5.4.1(b), the mixed LCP system for $d^{(2)}y$ dictates that for any index $i \in \gamma(\overline{z}) \cup d\beta_g(\overline{x};dx)$, the value is fixed to 0.
Since our single index $i=1$ is in $d\beta_g$, the formulation is simply:
$d^{(2)}y = 0$

(d) We solving for $y^{(2)}$.
Based on the formulation above, no matrix inversion or inequality solving is required. The second-order directional derivative is strictly determined by the index set classification:
$y^{(2)}(\overline{x}; 1) = 0$

\subsection*{Solutions for proof questions}

\paragraph{Question 1: AVI Constrained Mathematical Programs}
(a) The Mangasarian-Fromovitz Constraint Qualification (MFCQ) requires the existence of a direction that strictly satisfies the active inequality constraints while maintaining the equality constraints. In an MPAEC, the complementarity condition $\lambda^T(Dx+Ey+b)=0$, paired with $\lambda \ge 0$ and $Dx+Ey+b \le 0$, forces the feasible region to lie strictly on the boundary of the constraint set. Geometrically, this creates a tangent cone $\mathcal{T}(\overline{z}; \mathcal{F})$ that is a non-convex union of polyhedra, lacking the interiority required for the MFCQ to hold at feasible points like the origin.

(b) The nonlinear complementarity constraint $\lambda^T(Dx+Ey+b)=0$ can be reformulated using disjunctive linear constraints (i.e., either $\lambda_i = 0$ or $(Dx+Ey+b)_i = 0$). When the problem is viewed through a piecewise programming lens, the constraints are strictly linear within each active branch. Because the second derivatives of linear functions are identically zero, they do not contribute to the Hessian matrix of the nonlinear programming Lagrangean.

\paragraph{Question 2: NCP Constraints and Strong Regularity}

(a) The Strong Regularity Condition (SRC) at $\overline{z}$ requires two things:

The partial Jacobian matrix $\nabla_{\alpha(\overline{z})}F_{\alpha(\overline{z})}(\overline{z})$ is nonsingular.
The Schur complement of this matrix, defined as:
$$\nabla_{\beta(\overline{z})}F_{\beta(\overline{z})}(\overline{z}) - \nabla_{\alpha(\overline{z})}F_{\beta(\overline{z})}(\overline{z})(\nabla_{\alpha(\overline{z})}F_{\alpha(\overline{z})}(\overline{z}))^{-1}\nabla_{\beta(\overline{z})}F_{\alpha(\overline{z})}(\overline{z})$$
must be a P-matrix.

(b) This restriction is necessary in the standard multiplier-based approach because the derivation relies on the active set of constraints remaining stable along the critical direction $dz$. The strict complementarity restriction ensures that the index set $d\beta_b(\overline{x}; dx)$ is empty. Without it, the implicit function $y(x)$ might hit points of non-differentiability where the active constraint branch shifts. This combinatorial complication makes it impossible to guarantee that the inequality $dz^T(\nabla^2 f(\overline{z}) - \sum \pi_i \nabla^2 F_i(\overline{z}))dz \ge 0$ holds universally for all critical vectors under a single multiplier $\pi$.

\paragraph{Question 3: Implicit Programming Approach}

(a) The limit definition is:

$$y^{(2)}(\overline{x}; dx) = 2 \lim_{\tau \to 0+} \frac{y(\overline{x} + \tau dx) - \overline{y} - \tau y'(\overline{x}; dx)}{\tau^2}$$
(b) To prove positive homogeneity of degree 2, first note that for any scalar $\tau > 0$, the index sets scale such that $d\beta_s(\overline{x}; \tau dx) = d\beta_s(\overline{x}; dx)$ for $s \in \{a,b,g\}$. Furthermore, the quadratic form $dz^T \nabla^2 F_i(\overline{z}) dz$ is positively homogeneous of degree 2 in $dx$ (since $dz = (dx, y'(\overline{x};dx))$ and $y'$ is homogeneous of degree 1). Thus, scaling $dx$ by $\tau$ scales the constant vector in the mixed LCP system defining $y^{(2)}$ exactly by $\tau^2$. By the existence and uniqueness of the LCP solution, the solution itself scales by $\tau^2$. Finally, evaluating at 0 yields $y^{(2)}(\overline{x}; 0) = 0$.

(c) Proof: Use a Taylor expansion for $z(\tau) \equiv \overline{z} + \tau dz(\tau)$, where $dz(\tau) \equiv (dx, \frac{y(\tau)-\overline{y}}{\tau})$. Because $\overline{z}$ is a local minimum, we have for $\tau > 0$:

$$0 \le f(z(\tau)) - f(\overline{z}) = \nabla f(\overline{z})^T(z(\tau)-\overline{z}) + \frac{1}{2}(z(\tau)-\overline{z})^T \nabla^2 f(\overline{z})(z(\tau)-\overline{z}) + o(\tau^2)$$
Since $dz$ is a critical vector, the stationarity conditions imply $\nabla_x f(\overline{z})^T dx + \nabla_y f(\overline{z})^T y'(\overline{x};dx) = 0$. Subtracting this exact zero from the first-order term leaves $\nabla_y f(\overline{z})^T(dy(\tau) - y'(\overline{x};dx))$. Dividing the entire expansion by $\tau^2$ and taking the limit as $\tau \to 0+$ yields the final inequality directly via the definition of $y^{(2)}$.

Why it avoids multipliers: It completely circumvents $P(\overline{z})$ because it evaluates the objective function's curvature directly along the exact trajectory of the implicit function $y(x)$, rather than bounding the curvature using a Lagrangean dualized with KKT weights.

\paragraph{Question 4: Piecewise Programming and the Relaxed NLP}

(a) The relaxed NLP associated with $\overline{w}$ is:
\begin{equation*}
    \begin{array}{ll}
        \text{minimize }  &f(x,y) \\[4pt]
        \text{subject to }& Gx + Hy + a \le 0, L(x,y,\lambda) = 0,
    \end{array}
\end{equation*}
and the relaxed complementarity bounds are
\begin{equation*}
    \begin{cases}
        \lambda_i = 0, \quad g_i(x,y) \le 0 \quad \forall i \notin \mathcal{I}(\overline{z}),\\[4pt]
        \lambda_i \ge 0, \quad g_i(x,y) \le 0 \quad \forall i \in \mathcal{I}_0(\overline{z},\overline{\lambda}),\\[4pt]
        \lambda_i \ge 0, \quad g_i(x,y) = 0 \quad \forall i \in \mathcal{I}_+(\overline{z},\overline{\lambda}).
    \end{cases}
\end{equation*}
The relaxed critical cone $\mathcal{C}_{rl}(\overline{w})$ consists of all vectors $(dx, dy, d\lambda) \in \mathcal{T}(\overline{z}; Z) \times \mathbb{R}^l$ that are orthogonal to $\nabla f(\overline{z})$ and satisfy the linearized equations and inequalities of the active constraints defined above.

(b) Proof of uniqueness: Assume for contradiction there exists another multiplier $\lambda \in M(\overline{z}) \setminus \{\overline{\lambda}\}$. Because the multiplier set $M(\overline{z})$ is convex, we can choose a $\lambda$ arbitrarily close to $\overline{\lambda}$ such that $w = (\overline{z}, \lambda)$ falls inside the local neighborhood $V$.
Define the difference vector $dw = (0, 0, \lambda - \overline{\lambda})$. Because $z = \overline{z}$, the spatial difference $dz = (0,0)$.
Evaluating the quadratic form of the MPEC Lagrangean yields:
$$dw^T \nabla_{ww}^2\mathcal{L}^{MPEC}(\overline{w},\zeta,\pi,\eta)dw = dz^T \left( \nabla^2 f(\overline{z}) - \sum \pi_i \nabla^2 F_i(\overline{z}) \right) dz = 0$$
This evaluates to zero because $\mathcal{L}^{MPEC}$ contains no quadratic terms with respect to $\lambda$ (i.e., the block matrix $\nabla_{\lambda\lambda}^2 \mathcal{L}^{MPEC} \equiv 0$). However, if the Hessian is assumed to be strictly copositive on the critical cone, $dw^T \nabla_{ww}^2 \mathcal{L}^{MPEC} dw$ must be strictly greater than 0 for any non-zero $dw \in \mathcal{C}_{rl}(\overline{w})$. This contradiction proves that $d\lambda$ must be exactly 0, confirming that the multiplier $\overline{\lambda}$ is unique.

\section{Bibliographic remarks and Acknowledgment}

\textbf{This note is mainly based on Chapter 5, in the MPEC monograph of Zhi-Quan Luo, Jong-Shi Pang and Daniel Ralph.} Please see the relevant references: \cite{robinson2007local,robinson2009local,robinson2009local1,liu2025bidirectional,olikier2025projected,cui2026lipschitz,lin2023monotone,dussault2025b,wang2025analysis,huong2025generalized,robinson2009generalized,shapiro2021lectures,mordukhovich2023globally,besanccon2024flexible,wang2026algebraic,andreani2022optimality,burke1991exact,borghi2023constrained,wang2023sequential,strekalovsky2020local,benko2022sufficient,liang2025global,gu2020exact,mccormick1983nonlinear,fischetti2020branch,paradiso2025exact,robinson1980strongly,wang2025analysis1,shin2023near,shin2022exponential,dontchev2021lectures,bolte2024differentiating,wang2026damage,gfrerer2022local,nghia2025geometric,chen2025aubin,wang2023strong,pang1990newton,khanh2024globally,yu2026pattern,dussault2023exact,duy2023generalized,wang2023solving,shuo2026lecture,dussault2026polyhedral,liu1995sensitivity}. 

\bibliographystyle{unsrtnat}
\bibliography{reference5} 

\end{document}